\documentclass{amsart}
\usepackage{rotating}
\usepackage{enumerate}
\usepackage{amsfonts}
\usepackage{amsmath}
\usepackage{amssymb}
\usepackage{textcomp}
\usepackage{float,lscape}
\usepackage{mathtools}
\usepackage[titletoc,toc,title]{appendix}

\usepackage{tikz-cd}
\usepackage{enumitem}
\usepackage[arrow, matrix, curve, graph]{xy}
\usepackage{tikz}
\usetikzlibrary{arrows}
\usepackage[hyphens]{url}  
\newtheorem{theorem}{Theorem}[section]

\newtheorem{lemma}[theorem]{Lemma}

\newtheorem{corollary}[theorem]{Corollary}

\theoremstyle{definition}

\newtheorem{definition}[theorem]{Definition}

\newtheorem{example}[theorem]{Example}

\newtheorem{remark}[theorem]{Remark}

\newtheorem *{Theorem A}{Theorem A}
\newtheorem *{Theorem B}{Theorem B}
\newtheorem *{Corollary B}{Corollary B}
\newtheorem *{Corollary D}{Corollary D}

\newtheorem *{Problem1}{The group ring isomorphism problem [GRIP]}
\newtheorem *{Problem2}{The integral group ring isomorphism problem [IGRIP]}
\newtheorem *{Problem3}{Main Problem: The twisted group ring isomorphism problem [TGRIP]}

\newcommand{\RNum}[1]{\uppercase\expandafter{\romannumeral #1\relax}}
\newcommand{\rNum}[1]{\lowercase\expandafter{\romannumeral #1\relax}}

\numberwithin{equation}{section}

\begin{document}
	\title{Graded relations on crossed products}
\author{Ofir Schnabel}
\address{Department of Mathematics, ORT Braude College, 2161002 Karmiel, Israel}
\email{ofirsch@braude.ac.il}

	\begin{abstract}
		We classify crossed product gradings for arbitrary groups and fields up to several equivalence relations in terms of group actions and their orbits. 
	\end{abstract}

	\maketitle
	\bibliographystyle{abbrv}

	\section{Introduction}\label{Intro}\pagenumbering{arabic} \setcounter{page}{1}
	Throughout this note, $\mathbb{F}$ is a field, $G$ is a group and $A$ is an associative $\mathbb{F}$-algebra.
	Recall that a {\it grading} of $A$ by $G$ is a vector space decomposition
	\begin{equation}\label{eq:algebragrading}
		A=\bigoplus _{g\in G} A_g
	\end{equation}
	such that $A_gA_h\subseteq A_{gh}$. 
	
	Denote by $e$ the identity element in $G$. Then, it is straightforward that $A_e$ is an algebra. We call this sub-algebra of $A$, the \textit{base algebra} of the grading.
	
	The research of graded algebras received a lot of attention in the last decades (see e.g. \cite{EK13,NVO82}) and several equivalence relations were defined on graded algebras. 	A main goal was and still is to classify graded algebras up to these equivalence relations. 
	
	Probably the most common (see e.g. \cite{aljadeff2011simple,BSZ}) way to identify two gradings is by ``graded isomorphism" (see Definition~\ref{DefGradedIsomorphism}).	
	Roughly speaking, two $G$-graded algebras $A_G,B_G$ are graded isomorphic if there is an isomorphism between them which sends
	the homogeneous component $A_g$ to $B_g$ for every $g\in G$.
	
	However, other kinds of identification may sometimes be better off:
	\begin{itemize}
		\item When dealing with the intrinsic fundamental group $\pi_1(A)$ (see
		\cite[\S 3]{cibils2010} and \cite[\S 2.4]{ginosargradings}) of an
		algebra $A$ it is natural to study gradings of $A$ up to ``graded equivalence" (see Definition~\ref{DefGinosarEquivalent}). This relation is coarser than graded isomorphism (see also \cite{EK13}).  		
		\item In the study of structure (e.g. graded ideals,
		graded subspaces, radicals,\dots) or graded polynomial identities
		of graded algebras, the grading group itself does not play an
		important role, but can be replaced by any other group that
		realizes the same grading (see e.g. \cite{Clase, GordienkoSchnabel}). These gradings will be called ``weakly
		equivalent" (see \cite[Definition 1.4]{GordienkoSchnabel}).
	\end{itemize}
	
	In a sense, this paper is a continuation and generalization of \cite{karrer1973graded} where graded isomorphism classes of graded division algebras where described in terms of orbits under some group action. Here, we classify crossed products in similar terms. The group from \cite{karrer1973graded} is a subgroup of the group $K$ we are constructing and using here (see \S\ref{KgroupactingonGamma}). We were also motivated by the papers, \cite{BZ_real1, BZ_real2} where the authors classified the graded isomorphism classes and the graded equivalence classes of finite dimensional graded division algebras over the field of real numbers for abelian groups.
	
	In this paper, we put \cite{karrer1973graded} into a more general framework of crossed product gradings, as a particular case we classify graded division algebras up to graded isomorphism. As a consequence we classify graded division algebras over the field of real numbers up to graded isomorphism.
	
	A $G$-graded algebra~\eqref{eq:algebragrading} is a \textit{graded division algebra} if $A_g\neq \{0\}$ for any $g\in G$ and any $0\neq a \in A_g$ is an invertible element.
	A $G$-grading~\eqref{eq:algebragrading} is a \textit{crossed product} if $A_g$ admits an invertible element for any $g\in G$. Therefore, a $G$-graded division algebra is a particular case of a crossed product.
	
	For crossed products there are two more natural equivalences of gradings:
	\begin{itemize}
		\item Crossed products are used as ambient spaces in coding theory (see e.g. \cite{ginosarmoreno}), mainly when the grading group is finite. The ``alphabet" of a code consists of the elements of the component $A_e$.
		
		The appropriate identification in this regard is through ``graded isometry", which fixes the alphabet, but permits an isomorphism of the grading group (see Definition~\ref{def:Risometric}).
		\item In the context of realization-obstruction exact sequences (see e.g.\cite{CegarraGarzon, ginosar22}) the appropriate identification is through ``equivalent as Clifford systems" (see Definition~\ref{def:Cliffordequiv}) which is a graded isomorphism which additionally takes into account ``how" $A_e$ is embedded in $A$. This is a refinement of all the other mentioned equivalence relations.
	\end{itemize} 
	By its definition a crossed product is
	a free $R=A_e$ module with basis of invertible homogeneous elements $\{u_g\}_{g\in G}$, which  determines a map $\eta :G\rightarrow \text{Aut}_{\mathbb{F}}(R)$ by
	$$\eta (g)(r):=u_gru_g^{-1}$$
	for $r\in R$.
	It is easy to check that
	$$\eta (g)\eta (h)\eta (gh)^{-1}\in \text{Inn}(R)$$
	for any $g,h\in G$.
	Such a map $\eta:G\rightarrow $Aut$_{\mathbb{F}}(R)$ is called an \textit{outer action} of $G$ on $R$ and denote the set of all outer actions of $G$ on $R$ by Out$(G,R)$.
	The crossed product and the choice of a basis $\{u_g\}_{g\in G}$ also determines a map  
	$\alpha :G\times G\rightarrow R^*$ by
	$$\alpha(g,h):=u_gu_hu_{gh}^{-1}$$
	for any $g,h\in G$.
	The associativity of the crossed product yields 
	\begin{equation}\label{eq:etatwistcondition}
		\alpha (g_1,g_2)\alpha (g_1g_2,g_3)=\eta (g_1)(\alpha (g_2g_3))\alpha (g_1,g_2g_3)
	\end{equation}
	for any $g_1,g_2,g_3\in G$.
	Also, comparing 
	$u_gu_hru_h^{-1}u_g^{-1}$
	with $u_{gh}ru_{gh}^{-1}$
	yields 
	\begin{equation}\label{eq:etatwistcondition2}
		\eta (g)\circ \eta (h)=\alpha (g,h)\eta (gh)\alpha (g,h)^{-1}.
	\end{equation}
		
	With the above notation, for a fixed outer action $\eta \in \text{Out}(G,R)$ we will call a map $\alpha :G\times G\rightarrow R^*$ an $\eta$-twisting if~\eqref{eq:etatwistcondition} is satisfied. Alternatively, for a fixed $\alpha :G\times G\rightarrow R^*$ we will call a map $\eta :G\rightarrow \text{Aut}_{\mathbb{F}}(R)$ an $\alpha$-outer action if~\eqref{eq:etatwistcondition2} is satisfied.
	
	For a group $G$ and an algebra $R$, denote by $(R^*)^{G\times G}$ the set of functions from $G\times G$ to $R^*$ and define a subset $\Gamma$ of $(R^*)^{G\times G}\times \text{Out}(G,R)$ by
	$$\Gamma:=\Gamma_{G,R}=\{(\alpha,\eta)| \text{ s.t. }\alpha \text{ and } \eta  \text{ satisfying }~\eqref{eq:etatwistcondition} \text{ and }~\eqref{eq:etatwistcondition2}\}.$$
		For $(\alpha,\eta)\in \Gamma$ we denote the corresponding crossed product by $R^{\alpha}_{\eta}G$. We also notice that group rings, twisted group rings and skew group rings are particular cases of crossed products.
	\begin{example}\label{exintro}
		Let $G=C_2=\langle g\rangle$, $\mathbb{F}=\mathbb{R}$ and $R=\mathbb{C}$. $\text{Aut}_{\mathbb{R}}(\mathbb{C})$ admits two automorphisms, namely the trivial automorphism and the complex conjugation. With abuse of notation we will say either $\eta$ is trivial or $\eta$ is the complex conjugation. For both we assume  $$\alpha(1,1)=\alpha(g,1)=\alpha(1,g)=1.$$
		\begin{itemize}
			\item For $\alpha(g,g)=1$, which means in particular $u_g^2=1$, and $\eta$ trivial we get the group algebra $\mathbb{C}G=\mathbb{C}\oplus \mathbb{C}$.
			\item For $\alpha(g,g)=1$, which means in particular $u_g^2=1$, and $\eta$ is the complex conjugation we get
			$\mathbb{C}_{\eta}G\cong M_2(\mathbb{R})$.
			\item For $\alpha(g,g)=-1$, which means in particular $u_g^2=-1$, and $\eta$ is the complex conjugation we get
			$\mathbb{C}^{\alpha}_{\eta}G\cong\mathbb{H}$ the real Quaternions algebra. 
		\end{itemize} 
	\end{example}
	
	We notice that $\Gamma$ is ``too fine" in the sense that different basis choices to ``equivalent" crossed products give rise to different pairs in $\Gamma$. The objective of Theorem A is to give a criterion for distinct pairs in $\Gamma$ to correspond to ``equivalent" crossed products where here ``equivalent" depends on the equivalence relation we refer to. 
	
	A main tool for that is a group $K$ acting on $\Gamma$. We consider the group $(R^*)^G$ with pointwise multiplication, that is $(\lambda_1\cdot \lambda_2)(g)=\lambda_1 (g)\cdot \lambda_2 (g)$  for $\lambda _1,\lambda _2\in (R^*)^G$ and $g\in G$. 
	For $\lambda \in (R^*)^G, \varphi \in \text{Aut}_{\mathbb{F}}(R),\phi\in \text{Aut}(G))$ we
	denote
	$$\varphi(\lambda) (g):=\varphi(\lambda (g)) \text{ and } \phi(\lambda) (g):=\lambda (\phi ^{-1}(g)).$$
	We define in \S\ref{KgroupactingonGamma} a group
	$$K:=K_{G,R}=(R^*)^G\rtimes (\text{Aut}_{\mathbb{F}}(R)\times \text{Aut}(G))$$
	and furthermore we define an action~\eqref{eq:TheActionOfK} of $K$ on $\Gamma$ (see also Theorem~\ref{th:KisactingasagrouponGamma}).
	
	Then, with respect to the notion of equivalent as Clifford systems (see convention after Definition~\ref{def:Cliffordequiv}), graded isomorphism, graded isometric and graded equivalent, and with respect to the group $K$ above, we prove. 
	\begin{Theorem A}
		For an algebra $R$ and a group $G$, the group $K$ described in \S\ref{KgroupactingonGamma} is acting on $\Gamma$ by~\eqref{eq:TheActionOfK} such that two crossed products
		$R^{\alpha}_{\eta}G$ and $R^{\alpha^{\shortmid}}_{\eta ^{\shortmid}}G$ 
		\begin{enumerate}
			\item are equivalent as Clifford systems if and only if  $(\alpha,\eta)$ and $(\alpha ^{\shortmid},\eta ^{\shortmid})$ are in the same $(R^*)^G$ orbit.
			\item are graded isomorphic if and only if $(\alpha,\eta)$ and $(\alpha ^{\shortmid},\eta ^{\shortmid})$ are in the same $(R^*)^G\rtimes$ Aut$_{\mathbb{F}}(R)$ orbit. 
			\item are graded isometric if and only if $(\alpha,\eta)$ and $(\alpha ^{\shortmid},\eta ^{\shortmid})$ are in the same\\ $(R^*)^G\rtimes$ Aut$(G)$ orbit. 
			\item are graded equivalent if and only if $(\alpha,\eta)$ and $(\alpha ^{\shortmid},\eta ^{\shortmid})$ are in the same $K$ orbit. 
		\end{enumerate} 
	\end{Theorem A}
	By \cite[Corollary 3.6]{GordienkoSchnabel}, two crossed products are weakly graded equivalent (see \cite[Definition 1.4]{GordienkoSchnabel}) if and only if they are graded equivalent. Hence Theorem A above gives also the condition for crossed products to be weakly graded equivalent.
	
	We mention that for $R$ commutative, Theorem A $(3)$ appears in \cite[Theorem 3.5]{ginosarmoreno}. Also, Theorem A $(2)$  for the particular case of graded division algebras appears in \cite{karrer1973graded}.
	
	Denote by 
	$$\Psi=\Psi (G,R)=\{\eta \in \text{Out}(G,R)|\hspace{0.5cm} \exists \alpha \text{ such that } (\alpha ,\eta)\in \Gamma\}.$$
	It is important to mention that in general $\Psi$ is a proper subgroup of Out$(G,R)$. In fact there is a known criterion \cite[Theorem 3.4]{CegarraGarzon} for $\eta \in$Out$(G,R)$ to be in $\Psi$ (see also \cite[Theorem 7.1]{eilenberg1947cohomology}).
	For a fixed $\eta\in \Psi$, we denote
	$$\tilde{Z}^2_{\eta}(G,R^*)=\{\alpha \in (R^*)^{G\times G}\mid (\alpha, \eta)\in \Gamma\}.$$
	 Notice that if $R$ is commutative then this set is the group of 2-cocycles $Z^2_{\eta}(G,R^*)$. Then, we define an equivalence relation $\sim _{\eta}$ on $\tilde{Z}^2_{\eta}(G,R^*)$ by $\alpha \sim _{\eta} \alpha^{\shortmid}$ if $(\alpha, \eta)$ and $(\alpha^{\shortmid}, \eta)$ are in the same $(R^*)^G$-orbit (with respect to the action of $K$ on $\Gamma$), that is (by Theorem A) if $R^{\alpha}_{\eta}G$ and $R^{\alpha^{\shortmid}}_{\eta}G$ are equivalent as Clifford systems. Since the center $\mathcal{Z}(R)$ of $R$ is stable under any $\varphi\in $Aut$(R)$ we can restrict $\varphi $ to $\mathcal{Z}(R)$. Therefore we can also restrict an outer action $\eta\in \text{Out}(G,R)$ to $\mathcal{Z}(R)$. We will denote this restriction by $\bar{\eta}$ and notice that $\bar{\eta}$ is a homomorphism.
	Then, there is a well defined restriction of $\sim_{\eta}$ to $\sim _{\bar{\eta}}$ the induced relation on the second cohomology group $H^2_{\bar{\eta}}(G,\mathcal{Z}(R)^*)$. For $R$ commutative $\eta=\bar{\eta}$ and in those cases we abuse notation and write $H^2_{\eta}(G,R^*)$.
	
	Define an equivalence relation $\Delta$ on $\Psi$ as follows. For $\eta _1,\eta_2 \in \Psi$ we say $\eta_1 \Delta \eta_2$ if there exists $\varphi\in $Aut$(R)$ such that 
	$$\varphi\eta_1 (g)\varphi ^{-1}\eta_2(g)^{-1}\in \text{Inn}(R)$$ 
	for any $g\in G$.
	Denote also, for a fixed $\eta \in \Psi$, by $\Omega _{\eta,G}$
	the set of all crossed-products $R^{\alpha}_{\eta}G$ up to equivalence as Clifford systems. By \cite[Theorem 3.5]{CegarraGarzon} $H^2_{\bar{\eta}}(G,(\mathcal{Z}(R))^*)$ is acting freely and transitively on $\Omega _{\eta,G}$. It follows from the above and from Theorem A that there is a (non-canonical) identification between the set of all crossed-products $R^{\alpha}_{\eta}G$ up to graded isomorphism and $H^2_{\bar{\eta}}(G,\mathcal{Z}(R)^*)/\sim _{\bar{\eta}}$. 

	Let $G$ be a finite group. We would like to use the above to classify finite dimensional $G$-graded division algebras.
	It is straightforward to show that a graded algebra is a graded division algebra if and only if it is crossed product with $A_e=D$ a division algebra, so any $G$-graded division algebra can be written as $D^{\alpha}_{\eta}G$. So, using Theorem A we can classify in similar terms and with additional term of the Brauer group all graded division algebras up to the above mention equivalence relations.
	With the above setup, for a field $\mathbb{F}$, denote the base algebra of a graded division algebra by $D$ and denote the field $\mathcal{Z}(D)$ by $L$ and the Brauer group of $L$ by Br$(L)$. Then, the for a fixed $\mathbb{F}$ and $G$ the $G$-graded division algebras over $\mathbb{F}$ graded isomorphism classes can be classified by 
	\begin{equation}\label{eq:classify}
		\coprod_{[L:\mathbb{F}]<\infty} \coprod_{[D]\in \text{Br}(L)}\coprod_{[\eta]\in \Psi (D,G) /\Delta} \Omega _{\eta,G}.
	\end{equation}
	Lastly, we want to focus on real graded division algebras (see \cite{BZ_real1,BZ_real2,karrer1973graded,Lavit}), that is, with the notation of~\eqref{eq:classify}, $\mathbb{F}=\mathbb{R}$.
	Then $L$ is either 
	$\mathbb{R}$ or $\mathbb{C}$.
	\begin{enumerate}
		\item For $L=\mathbb{R}$, $D$ is either $\mathbb{R}$ or $\mathbb{H}$.
		For $D=\mathbb{R}$, $\eta$ is trivial, and for
		$D=\mathbb{H}$, by Skolem-Noether Theorem, $\eta$ is again trivial 
		\item For $L=\mathbb{C}$, $D$ is also $\mathbb{C}$. We have two automorphisms, the trivial automorphism and the complex conjugation. Hence, the kernel $N$ of $\eta$ is either $G$ itself or a subgroup of index $2$.
	\end{enumerate}
	Then, we have the following corollary.
	\begin{Corollary B}(see \cite[\S 8]{karrer1973graded} and \cite{Lavit})
		For a  graded $\mathbb{R}$-division algebra with a base algebra $R$ over a finite group $G$, $R$ is either 
		$\mathbb{R},\mathbb{H}$ or $\mathbb{C}$ and
		\begin{enumerate}
			\item For $R=\mathbb{R}$ the graded isomorphism classes are exactly the twisted group algebras $\mathbb{R}^{\alpha}G$ where $\alpha$ runs over the distinct classes in $H^2(G,\mathbb{R}^*)$.
			\item For $R=\mathbb{H}$, again there is a 1-1 correspondence between the graded isomorphism classes and the distinct classes in $H^2(G,\mathbb{R}^*)$.
			\item For $R=\mathbb{C}$, the kernel of $\eta$ is $N$ where $N$ is either $G$ or a subgroup of index $2$. In other words, there is an identification between the non-trivial actions and the  subgroups of $G$ of index $2$. For each $N$ there is a 1-1 correspondence between the graded isomorphism classes and the distinct sets $\{[\alpha],[\bar{\alpha}]\}$ in $H^2_{\eta}(G,\mathbb{C}^*)$.
			\end{enumerate}
	\end{Corollary B}

	In \S\ref{pre} we will introduce the definitions and give some background that will be needed in the sequel. 
	In \S\ref{actionofKonGamma} we define an action of $K$ on $\Gamma$ and show in Theorem~\ref{th:GammaisKinvariant} and Theorem~\ref{th:KisactingasagrouponGamma} that $K$ is acting as a group on $\Gamma$.
	In \S\ref{main} we prove Theorem A. 
	The last section is devoted to some interesting examples demonstrating the distinct equivalence relations and emphasizing the difference between them.
	\section{Relations on graded algebras}\label{pre}
	In this section we will recall the definitions of the equivalence relations mentioned in the introduction. 
	Let \begin{equation}\label{EqTwoGroupGradings}\Gamma_1 \colon A=\bigoplus_{g \in G} A_g,\qquad \Gamma_2
		\colon B=\bigoplus_{h \in H} B_h
	\end{equation} be two gradings  where $G$
	and $H$ are groups and $A$ and $B$ are $F$-algebras.
	
	\begin{definition}\label{DefGradedIsomorphism}
		The gradings~(\ref{EqTwoGroupGradings}) are \textit{isomorphic} if $G=H$ and there exists an isomorphism $\varphi \colon A \mathrel{\widetilde\to} B$
		of algebras such that $\varphi(A_g)=B_g$
		for all $g\in G$.
		In this case we say that $A$ and $B$ are \textit{graded isomorphic}.
	\end{definition}
	
	In some cases, such as in~\cite{ginosargradings}, less restrictive requirements are more
	suitable.
	\begin{definition}\label{DefGinosarEquivalent}
		The gradings~(\ref{EqTwoGroupGradings})  are \textit{equivalent} if there exists an isomorphism
		$\varphi \colon A \mathrel{\widetilde\to} B$
		of algebras and an isomorphism $\psi \colon G \mathrel{\widetilde\to} H$
		of groups such that $\varphi(A_g)=B_{\psi(g)}$
		for all $g\in G$.
	\end{definition}
	\begin{remark}
		The notion of graded equivalence was considered in~\cite[Remark after Definition 3]{BSZ}.
		In~\cite{Mazorchuk} it appears under the name of graded isomorphism and in~\cite[Section~3.1]{EK13} it appears under the name a weak isomorphism
		of gradings.
		More on differences in the terminology in
		graded algebras can be found in \cite[\S 2.7]{ginosargradings}.
	\end{remark}
	For two $G$-graded algebras~\eqref{EqTwoGroupGradings}, it is straightforward that a graded isomorphism or graded equivalence
	induce, by restriction, an algebra isomorphism between the base algebras $A_e,B_e$, and hence when considering graded isomorphism or graded equivalence between crossed-products we may assume that the grading group $G$ is the same and the base ring $R$ is the same.
	We are left with the problem of determining necessary and sufficient conditions for the existence of graded isomorphism or graded equivalence between $R^{\alpha _1}_{\eta _1}G$ and $R^{\alpha _2}_{\eta _2}G$, for a given group $G$, a given ring $R$, and $(\alpha_1,\eta _1),(\alpha_2,\eta _2) \in \Gamma$.

	Recall that a $G$ algebra grading~\eqref{eq:algebragrading} is a \textit{strong grading}
	if $A_gA_h=A_{gh}$ for any $g,h\in G$.
	\begin{definition}
		Let $G$ be a group and let $R$ be an algebra. A $G$-Clifford extension of $R$ is a pair $(A,\iota)$ such that $A$ is $G$-strongly graded and $\iota : R\hookrightarrow A$ is an embedding such that $\iota (R)=A_e$. 
	\end{definition}
	Notice that A $G$-Clifford extension takes into account, not only the graded structure of $A$, but also ``how" $A_e$ is embedded in $A$. So, it is natural that an equivalence of $G$-Clifford extensions will be a refinement of graded isomorphism.
	\begin{definition}\label{def:Cliffordequiv}
		Two $G$-Clifford system extensions $(A_1,\iota _1)$,  $(A_2,\iota _2)$ are equivalent if there exists a graded isomorphism $\psi :A_1\rightarrow A_2$ such that $\psi \circ \iota_1=\iota_2$.
	\end{definition}
	Clearly a $G$-crossed product is $G$-strongly graded.
	
	\underline{\textbf{convention}}: When addressing a crossed product $R^{\alpha}_{\eta }G$ with basis $\{u_g\}_{g\in G}$ as a $G$-Clifford extension of $R$ we mean with respect to the embedding $r\mapsto r\cdot 1$ for any $r\in R$. We have.
	\begin{corollary}
		Two crossed products 
		$R^{\alpha _1}_{\eta _1}G$ and $R^{\alpha _2}_{\eta _2}G$ are equivalent as Clifford systems if and only there exist a graded isomorphism $\psi :R^{\alpha _1}_{\eta _1}G\rightarrow R^{\alpha _2}_{\eta _2}G$
		such that the restriction of $\psi$
		to $R$ is trivial, that is $\psi (r)=r$ for all $r\in R$.
	\end{corollary}
	Similarly, we define
	\begin{definition}\label{def:Risometric}(see \cite{ginosarmoreno})
		Two crossed products 
		$R^{\alpha _1}_{\eta _1}G$ and $R^{\alpha _2}_{\eta _2}G$ are called $R$-isometric if there exists a graded equivalence $\psi :R^{\alpha _1}_{\eta _1}G\rightarrow R^{\alpha _2}_{\eta _2}G$
		such that the restriction of $\psi$
		to $R$ is trivial, that is $\psi (r)=r$ for all $r\in R$.
	\end{definition}

	We add the following diagram to help understand the order of coarsening between the different equivalence relations on crossed products.
	
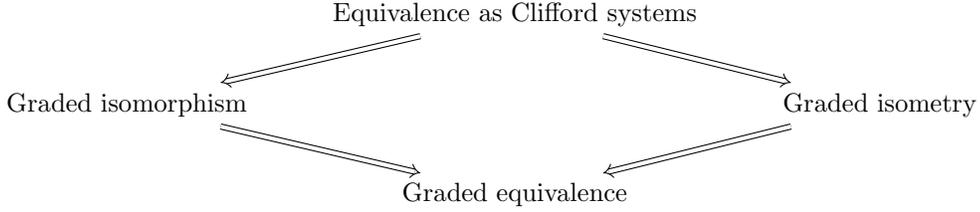
\begin{center}
		\begin{tikzcd}
			 & \text{Equivalence as Clifford systems}\arrow[dl, Rightarrow]\arrow[dr, Rightarrow]\\
			 \text{Graded isomorphism}\arrow[dr, Rightarrow] & & \text{Graded isometry}\arrow[dl, Rightarrow]\\
			 & \text{Graded equivalence}	
		\end{tikzcd}
	\end{center}
	\section{The group $K$ and its action on $\Gamma$}\label{actionofKonGamma}
	The main objective of this section is to construct the group $K$ mentioned in the introduction, which can be written as semi-direct product 
	$$K=(R^*)^G\rtimes (\text{Aut}_{\mathbb{F}}(R)\times \text{Aut}(G)),$$
	and to prove that $K$ is acting on $\Gamma$.
	It turns out that all the graded equivalence relations we define in the previous section can be identified with orbits of $\Gamma$ with respect to the restriction of the action of $K$ to certain subgroups of $K$.
	
	We break this section into two subsections. In the first we will give all the results that will be needed later and in the second part we will prove all these results.
	\subsection{$K$ is a group which acts on $\Gamma$}\label{KgroupactingonGamma}
	We start with the definition of $K$ by introducing an action of $\text{Aut}_{\mathbb{F}}(R)$ on $(R^*)^G$ and an action of $\text{Aut}(G))$ on 
	$(R^*)^G$.
	Let $\lambda \in (R^*)^G, \varphi \in \text{Aut}_{\mathbb{F}}(R),\phi\in \text{Aut}(G))$.
	Denote 
	$\lambda \mapsto \varphi (\lambda)$, $\lambda \mapsto \phi (\lambda)$ where

	$$\varphi(\lambda) (g):=\varphi(\lambda (g)) \text{ and } \phi(\lambda) (g):=\lambda (\phi ^{-1}(g)).$$
	\begin{lemma}\label{lemma:Kisagroup}
		Let $\lambda,\lambda _1,\lambda _2\in (R^*)^G$,
		let $\phi,\phi _1,\phi _2\in \text{Aut}(G))$ and let $\varphi,\varphi _1,\varphi_2 \in \text{Aut}_{\mathbb{F}}(R)$.
		Then, with the above notations 
		\begin{enumerate}
			\item $(\phi_1 \cdot \phi_2)(\lambda)=\phi_1 (\phi_2 (\lambda))$.
			\item $\phi (\lambda _1\cdot\lambda _2)=\phi (\lambda _1)\cdot \phi (\lambda _2)$.
			\item $(\varphi_1\cdot  \varphi_2)(\lambda)=\varphi_1 (\varphi_2 (\lambda))$.
			\item $\varphi (\lambda _1\cdot \lambda _2)=\varphi (\lambda _1)\cdot \varphi (\lambda _2)$.
			\item $\phi (\varphi (\lambda))=\varphi (\phi (\lambda))$.
		\end{enumerate}
	\end{lemma}

	\begin{corollary}\label{cor:definitionofK}
		With the above notations, 
		$(\text{Aut}_{\mathbb{F}}(R)\times \text{Aut}(G))$ acts on
		$(R^*)^G$ by
		$$(\phi ,\varphi)(\lambda):g\mapsto \varphi(\lambda (\phi^{-1}(g)))$$
		and therefore
		$$K=(R^*)^G\rtimes (\text{Aut}_{\mathbb{F}}(R)\times \text{Aut}(G))$$ is a group.
		The defining relations in $K$ are
		\begin{equation}
			\varphi \cdot \lambda=\varphi (\lambda)\cdot \varphi,\quad 
			\phi \cdot \lambda =\phi (\lambda)\cdot \phi,\quad \varphi \cdot \phi=\phi \cdot \varphi
		\end{equation}
		for any  
		$\lambda \in (R^*)^G$,
		$\phi \in \text{Aut}(G))$ and  $\varphi \in \text{Aut}_{\mathbb{F}}(R)$.
	\end{corollary}
	We will use the following notation.
	\begin{definition}
		Let $r\in R^*$. Then $\iota_r$ is an inner automorphism of $R$ defined by
		$$\iota_r(s):=rsr^{-1}$$
		for any $s\in R$. 	
	\end{definition} 
	Then, it is straightforward to show that
	\begin{lemma}\label{lemma:UsefullLemma}
		Let $r,r_1,r_2\in R^*$ and let $\varphi \in$Aut$(R)$. Then,
		\begin{enumerate}
			\item $\iota _{r_1}\circ \iota_{r_2}=\iota_{r_1r_2}$.
			\item $\varphi \circ \iota _r=\iota _{\varphi (r)}\circ \varphi$.
		\end{enumerate}
	\end{lemma}
	
	Next, we like to define an action of $K$ on the set $\Gamma $. We will do this in several steps.
	We will start by defining the following functions.
	\begin{definition}\label{def:actionofK}
		Let $(\alpha ,\eta)\in \left((R^*)^{G\times G}\times \text{Out}(G,R)\right)$, let 
		$\lambda \in (R^*)^G$, let $\varphi \in 
		\text{Aut}_{\mathbb{F}}(R)$ and let $\phi \in 
		\text{Aut}(G))$. Then, for any $g,g_1,g_2\in G$ and $r\in R$ define
		$$\varphi (\alpha ,\eta):=(\varphi (\alpha) ,\varphi (\eta)) \text{ where}$$
		\begin{enumerate}
			\item 
			$((\varphi (\eta))(g))(r)=\varphi (\eta (g)(\varphi ^{-1}(r))).$
			\item 
			$(\varphi (\alpha))(g_1,g_2)=\varphi (\alpha (g_1,g_2))$.
		\end{enumerate}	
		$$\phi (\alpha ,\eta):=(\phi (\alpha) ,\phi (\eta)) \text{ where}$$
		\begin{enumerate}\addtocounter{enumi}{2}
			\item 
			$((\phi (\eta))(g))(r)=(\eta ((\phi ^{-1}(g))))(r)$.
			\item 
			$(\phi (\alpha))(g_1,g_2)=\alpha (\phi ^{-1}(g_1),\phi ^{-1}(g_2))$
		\end{enumerate}	
		$$\lambda (\alpha ,\eta):=(\lambda _{\eta} (\alpha) ,\lambda (\eta)) \text{ where}$$
		\begin{enumerate}\addtocounter{enumi}{4}
			\item 
			$((\lambda (\eta))(g))(r)=\lambda (g)((\eta (g))(r))\lambda (g)^{-1}:=((\iota_{\lambda}\eta )(g))(r) $.
			\item 
			$(\lambda _{\eta}(\alpha))(g_1,g_2)=\lambda (g_1)(\eta (g_1)(\lambda (g_2)))\alpha (g_1,g_2)\lambda (g_1g_2)^{-1}$.
		\end{enumerate}
	\end{definition}
	\begin{corollary}\label{cor:concreteactionofK}
		With the above notation
		$$(\lambda \varphi \phi)(\alpha ,\eta)=
		(\lambda _{(\varphi \phi)(\eta)}((\varphi \phi)(\alpha)),(\lambda \varphi \phi)(\eta)$$
		where
		\begin{equation}\label{eq:concreteactiononeta}
			((\lambda \varphi \phi)(\eta ))(g)=\iota _{\lambda (g)}(\varphi \eta (\phi ^{-1}(g))\varphi ^{-1})
		\end{equation}
		and
		\begin{equation}\label{eq:concreteactionalpha}
			\lambda _{(\varphi \phi)(\eta)}((\varphi \phi)(\alpha))(g,h)=
			\lambda (g)(\varphi (\eta (\phi^{-1}(g)))\varphi^{-1}(\lambda (h)))\varphi (\alpha (\phi^{-1}(g),\phi ^{-1}(h)))\lambda (gh)^{-1}.
		\end{equation}
	\end{corollary}
	Recall that $\Gamma \subseteq  (R^*)^{G\times G}\times \text{Out}(G,R)$. The next theorem says that $\Gamma $ is invariant with respect to the maps in the above definition.  
	\begin{theorem}\label{th:GammaisKinvariant}
		With the above notations, for $(\alpha ,\eta)\in \Gamma$ we have that 
		\begin{enumerate}
			\item 
			$\varphi (\alpha ,\eta)\in \Gamma$.
			\item 
			$\phi (\alpha ,\eta)\in \Gamma$.
			\item 
			$\lambda (\alpha ,\eta)\in \Gamma$.
		\end{enumerate}
	\end{theorem}
	The following is immediate corollary.
	\begin{corollary}
		Let $(\alpha ,\eta)\in \Gamma$, let 
		$\lambda \in (R^*)^G$, let $\varphi \in 
		\text{Aut}_{\mathbb{F}}(R)$ and let $\phi \in 
		\text{Aut}(G))$.
		Then 
		$\lambda (\varphi (\phi (\alpha ,\eta)))\in \Gamma$.
	\end{corollary}
	By the above we can define the following map.
	\begin{definition}
		For any $\lambda \in (R^*)^G$, $\varphi \in 
		\text{Aut}_{\mathbb{F}}(R)$ and $\phi \in 
		\text{Aut}(G))$	
		we define a map
		$\psi :K\times \Gamma \to \Gamma$
		by 
		\begin{equation}\label{eq:TheActionOfK}
			\psi (\lambda \varphi \phi , (\alpha ,\eta))=
			\lambda (\varphi (\phi (\alpha ,\eta)))
		\end{equation}
		for any $(\alpha ,\eta)\in \Gamma$.
		We denote 
		$$
		\lambda \varphi \phi  (\alpha ,\eta):=\psi (\lambda \varphi \phi , (\alpha ,\eta))=\lambda (\varphi (\phi (\alpha ,\eta))).
		$$
	\end{definition}
	The next theorem says that the above is actually a group action of $K$ on $\Gamma $.
	
	\begin{theorem}\label{th:KisactingasagrouponGamma}
		For any $\lambda_1,\lambda_2 \in (R^*)^G$, $\varphi _1,\varphi_2\in 
		\text{Aut}_{\mathbb{F}}(R)$ and $\phi _1,\phi_2 \in 
		\text{Aut}(G))$	we have 
		$$
		(\lambda_1 \varphi_1 \phi_1)((\lambda_2 \varphi_2 \phi_2)  (\alpha ,\eta))=
		((\lambda_1 \varphi_1 \phi_1)(\lambda_2 \varphi_2 \phi_2))  (\alpha ,\eta))
		$$
		for any 
		$(\alpha ,\eta)\in \Gamma$.
	\end{theorem}
	\subsection{Proofs of the results from the previous section}
	
	\subsubsection{Proof of Lemma~\ref{lemma:Kisagroup}}
	\begin{enumerate}
		\item 
		$$
		\phi_1 (\phi_2 (\lambda))(g)= 
		\phi_2 (\lambda)(\phi _1^{-1}(g))=\lambda (
		\phi _2^{-1}(\phi _1^{-1}(g)))=\lambda ((\phi _1 \cdot \phi _2)^{-1}(g))=(\phi_1\cdot  \phi_2) (\lambda)(g).
		$$
		\item 
		$$(\phi (\lambda _1)\cdot \phi (\lambda _2))(g)=
		\lambda _1(\phi ^{-1}(g))\cdot \lambda _2(\phi ^{-1}(g))=(\lambda _1\cdot \lambda _2)(\phi ^{-1}(g))=(\phi (\lambda _1\cdot \lambda _2))(g).
		$$
		\item 
		$$\varphi_1 (\varphi_2 (\lambda))(g)=\varphi _1((\varphi_2(\lambda))(g))=\varphi _1(\varphi _2(\lambda (g)))=(\varphi_1 \cdot \varphi_2)(\lambda (g))=
		((\varphi_1 \cdot \varphi_2)(\lambda))(g).
		$$
		\item 
		$$
		\begin{aligned}
			& 
			(\varphi (\lambda _1 \cdot \lambda _2))(g)=
			\varphi ((\lambda _1 \cdot \lambda _2)(g))=\varphi 
			(\lambda _1(g)\cdot \lambda _2(g))=
			(\varphi 
			(\lambda _1))(g)\cdot (\varphi (\lambda _2))(g)=\\
			& =
			(\varphi 
			(\lambda _1)\cdot \varphi 
			(\lambda _2))(g)=
			(\varphi 
			(\lambda _1)\cdot \varphi 
			(\lambda _2))(g).
		\end{aligned}
		$$
		\item 
		$$(\phi (\varphi (\lambda)))(g)=\varphi (\lambda)(\phi ^{-1}(g))=\varphi (\lambda (\phi ^{-1}(g)))=\varphi (\phi (\lambda)(g))=(\varphi (\phi (\lambda)))(g).
		$$
	\end{enumerate}
	\qed 
	\subsubsection{Proof of Corollary~\ref{cor:concreteactionofK}}
	
	By Definition~\ref{def:actionofK} $(5)$
	$$
	\begin{aligned}
		& ((\lambda \varphi \phi)(\eta ))(g)=\lambda (g)(\varphi \phi (\eta)(g))\lambda ^{-1} (g)=\\
		& \iota _{\lambda (g)}\circ (\varphi \phi (\eta)(g))=
		\iota _{\lambda (g)}\circ ((\varphi \eta)(\phi^{-1}(g)))=\\
		& \iota _{\lambda (g)}\circ (\varphi \eta(\phi^{-1}(g))\varphi ^{-1})
	\end{aligned}
	$$
	proving~\eqref{eq:concreteactiononeta}.
	
	By the first part of the proof
	$$((\varphi \phi)(\eta))(g)=
	\varphi \eta (\phi ^{-1}(g))\varphi ^{-1}
	$$
	and hence by Definition~\ref{def:actionofK} $(6)$
	\begin{equation}
		\begin{aligned}
			& \lambda _{(\varphi \phi)(\eta)}((\varphi \phi)(\alpha))(g,h)=
			\lambda (g)((\varphi \phi)(\eta)(\lambda (h)))
			((\varphi \phi)(\alpha))(g,h)\lambda (gh)^{-1}=\\
			& \lambda (g)(\varphi (\eta (\phi^{-1}(g)))\varphi^{-1}(\lambda (h)))\varphi (\alpha (\phi^{-1}(g),\phi ^{-1}(h)))\lambda (gh)^{-1}
		\end{aligned}
	\end{equation}
	proving~\eqref{eq:concreteactionalpha}.
	\subsubsection{Proof of Theorem~\ref{th:GammaisKinvariant}}
	\begin{enumerate}
		\item 
		We need to show that the pair $(\varphi (\alpha) ,\varphi (\eta))\in \Gamma$, that is, to prove that $\varphi (\alpha)$ is a $\varphi (\eta)$-twisting (see~\eqref{eq:etatwistcondition}) and to prove that $\varphi (\eta)$ is a $\varphi (\alpha)$-outer action (see~\eqref{eq:etatwistcondition2}). We will first show that
		\begin{equation}\label{eq:varphietatwistcondition}
			\varphi (\alpha)(g_1,g_2)\varphi (\alpha) (g_1g_2,g_3)\stackrel{?}{=}\varphi (\eta) (g_1)(\varphi (\alpha) (g_2g_3))\varphi (\alpha) (g_1,g_2g_3).
		\end{equation}

		Since $(\alpha ,\eta)\in \Gamma$ we know that for any $g_1,g_2,g_3\in G$
		\begin{equation}
			\alpha (g_1,g_2)\alpha (g_1g_2,g_3)=\eta (g_1)(\alpha (g_2,g_3))\alpha (g_1,g_2g_3).
		\end{equation}
		Applying $\varphi $ on both sides of this equation yields
		$$\varphi(\alpha (g_1,g_2))\varphi(\alpha (g_1g_2,g_3))=\varphi(\eta (g_1)(\alpha (g_2g_3)))\varphi(\alpha (g_1,g_2g_3)).$$
		By plugging $\varphi ^{-1}\varphi$ after $\eta (g_1)$ on the right hand side of the equation we get that the equality in~\eqref{eq:varphietatwistcondition} holds.
		
		Next, we need to show that
		\begin{equation}\label{eq:varphietatwistcondition2}
			\varphi (\eta) (g_1)\circ \varphi(\eta) (g_2)\stackrel{?}{=}\iota_{\varphi (\alpha) (g_1,g_2)}\varphi (\eta) (g_1g_2).
		\end{equation}
		For any $r\in R$, by using~\eqref{eq:etatwistcondition2} the right hand side gives us
		$$
		\begin{aligned}
			& 
			\iota_{\varphi (\alpha) (g_1,g_2)}\varphi (\eta) (g_1g_2):r\mapsto 
			\iota_{\varphi (\alpha(g_1,g_2))}\varphi (\eta(g_1g_2)(\varphi ^{-1}(r)))=
			\\
			&
			\iota_{\varphi (\alpha(g_1,g_2))}\varphi(\iota^{-1}_{\alpha(g_1,g_2)}\eta(g_1)\eta (g_2)(\varphi^{-1}(r)))=\ldots
		\end{aligned}
		$$ 
		and by using Lemma~\ref{lemma:UsefullLemma} this is equal to 
		$$\ldots=(\varphi (\eta) (g_1)\circ \varphi(\eta) (g_2))(r)$$
		which is the left hand side of~\eqref{eq:varphietatwistcondition2}.
		Thus, $(\varphi (\alpha),\varphi(\eta))\in \Gamma $. This concludes the proof of Theorem~\ref{th:GammaisKinvariant}(1).
		
		\item 
		Denote 
		$(\alpha^{\shortmid},\eta^{\shortmid})=\phi(\alpha,\eta)$ and recall that here
		$\alpha^{\shortmid}(g_1,g_2)=\alpha (\phi ^{-1}(g_1),\phi ^{-1}(g_2))$ and $\eta ^{\shortmid}(g)=\eta(\phi ^{-1}(g_1))$.
		Since $(\alpha, \eta)\in \Gamma$ it satisfies~\eqref{eq:etatwistcondition} for any 3 group elements, in particular for
		$\phi ^{-1}(g_1),\phi ^{-1}(g_2),\phi ^{-1}(g_3)\in G$ and therefore
		$$
		\begin{aligned}
			& \alpha^{\shortmid} (g_1,g_2)\alpha^{\shortmid} (g_1g_2,g_3)=\alpha (\phi ^{-1}(g_1),\phi ^{-1}(g_2))\alpha (\phi ^{-1}(g_1)\phi ^{-1}(g_2),\phi ^{-1}(g_3))=\\
			&=\eta (\phi ^{-1}(g_1))(\alpha (\phi ^{-1}(g_2)\phi ^{-1}(g_3)))\alpha (\phi ^{-1}(g_1),\phi ^{-1}(g_2)\phi ^{-1}(g_3))=\\
			& =
			\eta^{\shortmid} (g_1)(\alpha^{\shortmid} (g_2,g_3))\alpha^{\shortmid} (g_1,g_2g_3)
		\end{aligned}
		$$
		and hence $\alpha^{\shortmid}$ is an $\eta^{\shortmid}$-twisting (see~\eqref{eq:etatwistcondition}).

		Next, since $(\alpha,\eta)\in \Gamma$ and since $\phi^{-1}(g_1),\phi^{-1}(g_2)\in G$ we know that

		\begin{equation}
			\begin{aligned}
				& \eta^{\shortmid}(g_1)\circ \eta^{\shortmid}(g_2)=
				\eta(\phi^{-1}(g_1))\circ \eta(\phi^{-1}(g_2))=\\
				& =\alpha(\phi^{-1}(g_1),\phi^{-1}(g_2))\eta(\phi^{-1}(g_1)\phi^{-1}(g_2))\alpha(\phi^{-1}(g_1),\phi^{-1}(g_2))^{-1}=\\
				& =\iota_{\alpha^{\shortmid}(g_1,g_2)}\eta^{\shortmid}(g_1g_2).
			\end{aligned}
		\end{equation}
		Hence, $\eta^{\shortmid}$ is an  $\alpha^{\shortmid}$-outer action (see~\eqref{eq:etatwistcondition2}) and therefore
		$(\alpha^{\shortmid},\eta^{\shortmid})\in \Gamma$. This concludes the proof of Theorem~\ref{th:GammaisKinvariant}(2).
		\item 
		Denote 
		$(\lambda (\alpha),\lambda (\eta))=\lambda(\alpha,\eta)$.
		We want to show first that
		\begin{equation}
			(\lambda (\alpha)) (g_1,g_2)(\lambda (\alpha)) (g_1g_2,g_3)\stackrel{?}{=}(\lambda (\eta)) (g_1)((\lambda (\alpha))) (g_2,g_3))(\lambda (\alpha)) (g_1,g_2g_3).
		\end{equation}
		The left hand side of the equation is
		$$
		\begin{aligned}
			& 
			\lambda (g_1)\eta(g_1)(\lambda(g_2))\alpha (g_1,g_2)\lambda (g_1g_2)^{-1}\lambda (g_1g_2)\eta(g_1g_2)(\lambda (g_3))\alpha (g_1g_2,g_3)\lambda (g_1g_2g_3)^{-1}
			=\\
			& =
			\lambda (g_1)\eta(g_1)(\lambda(g_2))\alpha (g_1,g_2)\eta(g_1g_2)(\lambda (g_3))\alpha (g_1g_2,g_3)\lambda (g_1g_2g_3)^{-1}.
		\end{aligned}
		$$
		The right hand side of the equation is
		$$
		\begin{aligned}
			& 
			\lambda (g_1)\eta(g_1)(\lambda(g_2)\eta (g_2)(\lambda (g_3))\alpha (g_2,g_3)\lambda (g_2g_3)^{-1})\lambda (g_1)^{-1}\lambda (g_1)\eta (g_1)(\lambda (g_2g_3))\alpha (g_1,g_2g_3)\lambda(g_1g_2g_3)^{-1}=\\
			&
			=
			\lambda (g_1)\eta(g_1)(\lambda(g_2)\eta (g_2)(\lambda (g_3))\alpha (g_2,g_3)\lambda (g_2g_3)^{-1})\eta (g_1)(\lambda (g_2g_3))\alpha (g_1,g_2g_3)\lambda(g_1g_2g_3)^{-1}.
		\end{aligned}
		$$
		We see that we can eliminate from the left on both sides of the equation
		$\lambda (g_1)\eta(g_1)(\lambda(g_2))$
		and from the right we can eliminate
		$\lambda(g_1g_2g_3)^{-1}$.
		Then, after these eliminations, the left hand side is
		$$\alpha (g_1,g_2)\eta(g_1g_2)(\lambda (g_3))\alpha (g_1g_2,g_3)$$
		and the right hand side is
		$$
		\begin{aligned}
			& 
			\eta(g_1)(\eta (g_2)(\lambda (g_3))\alpha (g_2,g_3)\lambda (g_2g_3)^{-1})\eta (g_1)(\lambda (g_2g_3))\alpha (g_1,g_2g_3)=\\
			&=\eta(g_1)(\eta (g_2)(\lambda (g_3))\alpha (g_2,g_3))\alpha (g_1,g_2g_3)=\\
			&=\eta(g_1)(\eta (g_2)(\lambda (g_3)))\cdot \eta(g_1)(\alpha (g_2,g_3))\alpha (g_1,g_2g_3)
		\end{aligned}
		$$
		which is equal by~\eqref{eq:etatwistcondition2} to
		$$\alpha (g_1,g_2)\eta(g_1g_2)(\lambda (g_3))\alpha (g_1,g_2)^{-1}\eta (g_1)(\alpha(g_2,g_3))\alpha (g_1,g_2g_3)$$
		and again we can eliminate from the left 
		$\alpha (g_1,g_2)\eta(g_1g_2)(\lambda (g_3))$.
		So now we need to show that
		$$\alpha (g_1g_2,g_3)\stackrel{?}{=}\alpha (g_1,g_2)^{-1}\eta (g_1)(\alpha(g_2,g_3))\alpha (g_1,g_2g_3).$$
		Multiplying both from the left by $\alpha (g_1,g_2)$
		we get on the left side 
		$\alpha (g_1,g_2)\alpha (g_1g_2,g_3)$
		and on the right we get 
		$\eta (g_1)(\alpha(g_2,g_3))\alpha (g_1,g_2g_3)$.
		Since $(\alpha ,\eta)\in \Gamma$ these two expressions are equal by~\eqref{eq:etatwistcondition}.
		
		Next, we need to show that
		\begin{equation}\label{eq:NeedToProve}
			(\lambda (\eta))(g_1)\circ (\lambda (\eta))(g_2)\stackrel{?}{=}\iota_{(\lambda (\alpha))(g_1,g_2)}(\lambda (\eta))(g_1g_2).
		\end{equation}
		By the definition of $\lambda (\eta)$ 
		\begin{equation}\label{eq:firststep}
			\iota_{(\lambda (\alpha))(g_1,g_2)}(\lambda (\eta))(g_1g_2)=\iota_{(\lambda (\alpha))(g_1,g_2)}\iota_{\lambda(g_1g_2)}\eta(g_1g_2)=\ldots
		\end{equation}
		Since $\eta$ is an $\alpha $-outer action (see~\eqref{eq:etatwistcondition2}) this is equal to
		$$
		\ldots =\iota_{(\lambda (\alpha))(g_1,g_2)}\iota_{\lambda(g_1g_2)}\iota_{\alpha^{-1}(g_1,g_2)}\eta(g_1)\eta(g_2)=\iota_{[(\lambda (\alpha))(g_1,g_2)\lambda(g_1g_2)\alpha^{-1}(g_1,g_2)]}\eta(g_1)\eta(g_2)=\ldots
		$$ 
		which is equal, by the definition of $\lambda (\alpha)$ to 
		$$
		\begin{aligned}
			& \ldots =
			\iota_{[\lambda (g_1)\eta(g_1)(\lambda (g_2))\alpha(g_1,g_2)\lambda(g_1g_2)^{-1}\lambda(g_1g_2)\alpha^{-1}(g_1,g_2)]}\eta(g_1)\eta(g_2)=\\
			& \iota_{[\lambda (g_1)\eta(g_1)(\lambda (g_2))]}\eta(g_1)\eta(g_2)=
			\iota_{\lambda (g_1)}\iota _{\eta(g_1)(\lambda (g_2))}\eta(g_1)\eta(g_2)=\ldots
		\end{aligned}
		$$
		Then, by Lemma~\ref{lemma:UsefullLemma} this is equal to 
		\begin{equation}\label{eq:finalstep}
			\iota_{\lambda(g_1)}\eta(g_1)\iota_{\lambda(g_2)}\eta(g_2)=\lambda (\eta)(g_1)\circ \lambda (\eta)(g_2).
		\end{equation}
		To conclude, we proved equality between the left hand side of~\eqref{eq:firststep} and the right hand side of~\eqref{eq:finalstep}, that is we proved that the equality in~\eqref{eq:NeedToProve} holds.
		
		Consequently $(\lambda (\alpha),\lambda(\eta))\in \Gamma $ and this conclude the proof of Theorem~\ref{th:GammaisKinvariant}(3).
	\end{enumerate}
	\qed 
	
	\subsubsection{Proof of Theorem~\ref{th:KisactingasagrouponGamma}}
	The proof will be in several steps.
	Let  $\lambda ,\lambda_1,\lambda_2 \in (R^*)^G$, $\varphi, \varphi _1,\varphi_2\in 
	\text{Aut}_{\mathbb{F}}(R)$ and $\phi,\phi _1,\phi_2 \in 
	\text{Aut}(G))$	and let $(\alpha ,\eta)\in \Gamma$.
	We will prove that
	\begin{enumerate}
		\item  
		$(\lambda _1\lambda _2)(\alpha, \eta)=\lambda _1	(\lambda _2(\alpha, \eta))$.
		\item 
		$(\varphi _1\varphi _2)(\alpha, \eta)=\varphi _1	(\varphi _2(\alpha, \eta))$.
		\item 
		$(\phi _1\phi _2)(\alpha, \eta)=\phi _1	(\phi _2(\alpha, \eta))$.
		\item 
		$(\phi \circ \varphi )(\alpha, \eta)=(\varphi \circ \phi )(\alpha, \eta)$.
		\item 
		$\varphi  (\lambda (\alpha, \eta))=(\lambda \varphi )(\alpha, \eta)=\varphi (\lambda)(\varphi (\alpha, \eta))$.
		\item 
		$\phi  (\lambda (\alpha, \eta))=(\lambda \phi )(\alpha, \eta)=\phi (\lambda)(\phi (\alpha, \eta))$.
		\item 
		$(\lambda_1 \varphi_1 \phi_1)((\lambda_2 \varphi_2 \phi_2)  (\alpha ,\eta))=
		((\lambda_1 \varphi_1 \phi_1)(\lambda_2 \varphi_2 \phi_2))  (\alpha ,\eta))
		$.
	\end{enumerate}
	\textbf{Proof of step (1)}
	Denote $(\alpha^{\shortmid}, \eta^{\shortmid })=(\lambda _2(\alpha, \eta))$ and 
	$(\alpha^{\shortmid \shortmid}, \eta^{\shortmid \shortmid })=\lambda _1	(\alpha^{\shortmid}, \eta ^{\shortmid})=\lambda _1	(\lambda _2(\alpha, \eta))$.
	
	By the definition of $\lambda (\eta)$ we have that
	$$(\lambda _1\lambda _2) (g)(\eta):r\mapsto \iota _{\lambda _1\lambda _2 (g)}(\eta (g)(r))$$
	for 
	$g\in G$ and $r\in R$ while 
	$$\lambda _1	(\lambda _2(\eta))(g):r\mapsto \iota _{\lambda _1 (g)}(\eta ^{\shortmid} (g)(r))=
	\iota _{\lambda _1 (g)}( \iota _{\lambda _2 (g)}(\eta  (g)(r)))$$ and 
	hence $(\lambda _1\lambda _2) (\eta)=\lambda _1(\lambda _2 (\eta))$.
	
	Next, 
	\begin{equation}
		\begin{aligned}
			& \lambda_1(\lambda_2(\alpha (g_1,g_2)))=\\
			& =\alpha ^{\shortmid \shortmid}(g_1,g_2)=\lambda _1 (\alpha ^{\shortmid})(g_1,g_2)=
			\lambda _1(g_1)\eta ^{\shortmid}(g_1)(\lambda _1(g_2))\alpha ^{\shortmid}(g_1,g_2)\lambda _1^{-1}(g_1g_2)=\\
			& =\lambda _1(g_1)\iota_{\lambda _2(g_1)}(\eta (g_1)(\lambda _1(g_2)))\lambda _2(g_1)\eta (g_1)(\lambda _2(g_2))\alpha (g_1,g_2)\lambda _2^{-1}(g_1g_2)\lambda _1^{-1}(g_1g_2)=\\
			&=\lambda _1(g_1)\lambda _2(g_1)(\eta (g_1)(\lambda _1(g_2))\eta (g_1)(\lambda _2(g_2))\alpha (g_1,g_2)\lambda _2^{-1}(g_1g_2)\lambda _1^{-1}(g_1g_2)=\\
			&=(\lambda _1\lambda _2)(g_1)\eta (g_1)(\lambda _1\lambda _2(g_2))\alpha (g_1,g_2)(\lambda _1\lambda _2)(g_1g_2)^{-1}=\\
			&=(\lambda _1\lambda _2)(\alpha (g_1,g_2)).
		\end{aligned}
	\end{equation}
	This conclude the proof of step (1).

	\textbf{Proof of step (2)}
	Denote 
	$$(\alpha ^{\shortmid },\eta ^{\shortmid}):=\varphi_2 (\alpha,\eta) \text { and }      (\alpha ^{\shortmid \shortmid},\eta ^{\shortmid \shortmid}):=\varphi _1(\alpha ^{\shortmid },\eta ^{\shortmid})=\varphi_1(\varphi_2 (\alpha,\eta)).
	$$
	Then, 
	$$\alpha ^{\shortmid \shortmid}(g_1,g_2)=\varphi_1 (\alpha ^{\shortmid })(g_1,g_2)=\varphi_1 (\varphi_2 (\alpha (g_1,g_2)))=\varphi_1\varphi_2(\alpha (g_1,g_2)).$$
	Next,
	$$\eta ^{\shortmid \shortmid}(g)(r)=\varphi_1 (\eta ^{\shortmid}(g)(\varphi_1^{-1}(r)))=\varphi_1(\varphi_2(\eta (g)(\varphi_2^{-1}(\varphi_1^{-1}(r)))))=\varphi_1\varphi_2\eta (g)(\varphi_1\varphi_2)^{-1}(r).$$
	This conclude the proof of step (2).

	\textbf{Proof of step (3)}
	This is clear by the definition 
	$$(\phi (\alpha))(g_1,g_2)=\alpha (\phi ^{-1}(g_1),\phi ^{-1}(g_2)), \text{ and } (\phi (\eta))(g)=\eta (\phi ^{-1}(g)).$$
	
	\textbf{Proof of step (4)}
	Denote 
	$$ (\alpha^{\shortmid},\eta ^{\shortmid}):=\varphi (\alpha ,\eta), \text{ and }(\tilde{\alpha},\tilde{\eta}):=\phi (\alpha ,\eta).$$
	Then,
	$$((\phi \circ \varphi) (\alpha))(g_1,g_2)=
	(\phi (\alpha^{\shortmid}))(g_1,g_2)=\alpha ^{\shortmid}(\phi ^{-1}(g_1),\phi^{-1}(g_2))=\varphi ( \alpha (\phi ^{-1}(g_1),\phi^{-1}(g_2)))$$
	and
	$$((\varphi \circ \phi) (\alpha))(g_1,g_2)=
	(\varphi (\tilde{\alpha}))(g_1,g_2)=
	\varphi (\tilde{\alpha}(g_1,g_2))
	=\varphi ( \alpha (\phi ^{-1}(g_1),\phi^{-1}(g_2)))$$
	and so $(\phi \circ \varphi) (\alpha)=(\varphi \circ \phi) (\alpha)$.
	
	Now, 
	$$((\phi \circ \varphi) (\eta (g)))(r)=(\phi (\eta ^{\shortmid}(g)))(r)=(\eta ^{\shortmid}(\phi^ {-1}(g)))(r)=\varphi (\eta (\phi ^{-1}(g))(\varphi ^{-1}(r)))$$
	and
	$$((\varphi \circ \phi) (\eta (g)))(r)=
	\varphi (\tilde{\eta}(g)(\varphi ^{-1}(r)))
	=\varphi (\eta (\phi ^{-1}(g))(\varphi ^{-1}(r)))$$
	and so $(\phi \circ \varphi) (\eta)=(\varphi \circ \phi) (\eta)$.
	This conclude the proof of step (4).
	
	\textbf{Proof of step (5)}
	Denote 
	$$ (\alpha^{\shortmid},\eta ^{\shortmid}):=\lambda (\alpha ,\eta), \text{ and }(\tilde{\alpha},\tilde{\eta}):=\varphi (\alpha ,\eta).$$
	Then,
	$$((\varphi \lambda)(\alpha))(g_1,g_2)=
	\varphi (\alpha ^{\shortmid}(g_1,g_2))=\varphi (\lambda (g_1)\eta (g_1)(\lambda (g_2))\alpha (g_1,g_2)\lambda (g_1g_2)^{-1})$$
	and on the other hand by the definition of the multiplication in $K$ 
	$$
	\begin{aligned}
		& ((\lambda \varphi)(\alpha))(g_1,g_2)=
		\varphi (\lambda)(\varphi(\alpha))(g_1,g_2)=
		\varphi (\lambda)(\tilde{\alpha}(g_1,g_2))=\\
		& =\varphi (\lambda)(g_1)\tilde{\eta}(\varphi (\lambda)(g_2))\tilde{\alpha}(g_1,g_2)\varphi (\lambda)(g_1g_2)^{-1}=\\
		& = \varphi (\lambda)(g_1)\varphi (\eta (g_1)(\varphi ^{-1}(\varphi (\lambda)(g_2))))\varphi (\alpha (g_1,g_2))\varphi (\lambda (g_1g_2))^{-1}=\\
		& = \varphi (\lambda)(g_1)\varphi (\eta (g_1)(\lambda(g_2)))\varphi (\alpha (g_1,g_2))\varphi (\lambda (g_1g_2))^{-1}=\\
		& = \varphi (\lambda (g_1)\eta (g_1)(\lambda (g_2))\alpha (g_1,g_2)\lambda (g_1g_2)^{-1}).
	\end{aligned}
	$$
	We conclude that 
	$((\varphi \lambda)(\alpha))(g_1,g_2)=((\lambda \varphi)(\alpha))(g_1,g_2)$.
	
	Now, 
	$$(\varphi (\eta ^{\shortmid})(g):r\mapsto \varphi (\eta^{\shortmid}(g_1)(\varphi ^{-1}(r)))=\varphi (\lambda (g_1)\eta (g_1)\varphi ^{-1}(r)\lambda (g_1)^{-1})$$
	and on the other hand
	$$(\varphi (\lambda)\tilde{\eta})(g):r\mapsto 
	(\varphi (\lambda))(g)\tilde{\eta}(g)(r) (\lambda))(g))^{-1}=\varphi (\lambda (g))\varphi (\eta (g)(\varphi^{-1}(r)))\varphi (\lambda (g))^{-1}$$
	and consequently
	$((\varphi \lambda)(\eta))(g)=((\lambda \varphi)(\eta)(g)$.
	
	This conclude the proof of step (5).

	\textbf{Proof of step (6)}
	Denote 
	$$ (\alpha^{\shortmid},\eta ^{\shortmid}):=\lambda (\alpha ,\eta), \text{ and }(\tilde{\alpha},\tilde{\eta}):=\phi (\alpha ,\eta).$$
	Denote also $\tilde{\lambda}=\phi (\lambda)$.
	
	Then,
	$$
	\begin{aligned}
		&((\phi \lambda) \alpha )(g_1,g_2)= (\phi (\alpha ^{\shortmid})(g_1,g_2)=\alpha ^{\shortmid}(\phi ^{-1}(g_1),\phi ^{-1}(g_2))=\\
		& =
		\lambda (\phi ^{-1})\eta(\phi ^{-1}(g_1))(\lambda (\phi ^{-1}(g_2)))\alpha (\phi ^{-1}(g_1),\phi ^{-1}(g_2))\lambda (\phi ^{-1}(g_1)\phi ^{-1}(g_2))^{-1}
	\end{aligned}
	$$
	and on the other hand 
	$$
	\begin{aligned}
		& (( \lambda \phi) \alpha ^{\shortmid})(g_1,g_2)=(\phi (\lambda))(\phi (\alpha (g_1,g_2)))=(\tilde{\lambda}(\tilde{\alpha}))(g_1,g_2)=\\
		& \tilde{\lambda}(g_1)\tilde{\eta}(g_1)(\tilde{\lambda}(g_2))\tilde{\alpha}(g_1,g_2)\tilde{\lambda}(g_1g_2)^{-1}=\\
		& =
		\lambda (\phi ^{-1})\eta(\phi ^{-1}(g_1))(\lambda (\phi ^{-1}(g_2)))\alpha (\phi ^{-1}(g_1),\phi ^{-1}(g_2))\lambda (\phi ^{-1}(g_1)\phi ^{-1}(g_2))^{-1}
	\end{aligned}
	$$
	and we conclude that $((\phi \lambda) \alpha )(g_1,g_2)=(( \lambda \phi) \alpha )(g_1,g_2)$.
	
	Next, 
	$$((\phi \lambda)(\eta)(g)=\phi (\eta ^{\shortmid})(g):r\mapsto \eta ^{\shortmid}(\phi ^{-1}(g)(r))=\lambda (\phi^{-1}(g))\eta(\phi ^{-1}(g))(r)\lambda (\phi^{-1}(g))^{-1}$$
	and on the other hand
	$$
	\begin{aligned}
		&
		((\lambda \phi )(\eta)(g)=
		((\phi (\lambda))( \phi (\eta))(g)=
		(\tilde{\lambda}\tilde{\eta}) (g):r\mapsto 
		\tilde{\lambda}(g)\tilde{\eta}(g)(r)\tilde{\lambda}(g)^{-1}=\\
		& =\lambda (\phi^{-1}(g))\eta(\phi ^{-1}(g))(r)\lambda (\phi^{-1}(g))^{-1}
	\end{aligned}
	$$
	and so 
	$((\phi \lambda)(\eta)(g)=((\lambda \phi )(\eta)(g)$.
	
	This conclude the proof of step (6).
	
	\textbf{Proof of step (7)}
	By the definition of the multiplication in $K$
	$$
	\begin{aligned}
		&((\lambda _1\varphi _1\phi_1)(\lambda _2\varphi _2\phi_2))(\alpha ,\eta )=(\lambda _1\varphi _1\phi _1(\lambda _2)\varphi _1\varphi_2\phi_1\phi_2)(\alpha ,\eta )=\\
		&=
		\lambda _1\varphi _1\phi _1(\lambda _2)
		(\varphi_1\varphi_2(\phi_1\phi_2(\alpha ,\eta )))
	\end{aligned}
	$$
	and on the other hand 
	$$
	\begin{aligned}
		& (\lambda _1\varphi _1\phi_1)((\lambda _2\varphi _2\phi_2)(\alpha ,\eta ))=
		\lambda _1(\varphi _1(\phi_1(\lambda _2(\varphi _2(\phi_2(\alpha ,\eta ))))))=\\
		& =
		\lambda _1(\varphi _1(\phi_1(\lambda _2)(\phi_1(\varphi _2(\phi_2(\alpha ,\eta ))))))=
		\lambda _1(\varphi _1(\phi_1(\lambda _2))(\varphi_1(\phi_1(\varphi _2(\phi_2(\alpha ,\eta ))))))=\\
		&=\lambda_1\varphi_1\phi_1(\lambda_2)(\varphi_1(\varphi_2(\phi_1 (\phi_2 (\alpha,\eta)))))=\lambda _1\varphi _1\phi _1(\lambda _2)
		(\varphi_1\varphi_2(\phi_1\phi_2(\alpha ,\eta ))).
	\end{aligned}
	$$
	This conclude the proof of step (7) and the proof Theorem~\ref{th:KisactingasagrouponGamma}.
	\qed

	\section{Proof of Theorem A}\label{main}
	Let $\{u_g\}_{g\in G}$ be a basis for $R^{\alpha _1}_{\eta _1}G$ and let $\{v_g\}_{g\in G}$ be a basis for $R^{\alpha_2}_{\eta_2}G$.
	Assume $R^{\alpha_1}_{\eta_1}G$ and $R^{\alpha_2}_{\eta_2}G$ are graded equivalent. Then there is an isomorphism $\psi :R^{\alpha_1}_{\eta_1}G\rightarrow 	R^{\alpha_2}_{\eta_2}G$, $\phi \in $Aut$(G)$ and $\lambda \in (R^*)^G$ such that
	\begin{equation}\label{eq:thedefinitionofPsi}
		\psi \left(\sum_{g\in G}r_gu_g\right)=\sum_{g\in G}\varphi(r_g)\lambda(\phi (g))^{-1}v_{\phi (g)},
	\end{equation}
	where here $\varphi \in $Aut$(R)$ is the restriction of $\psi$ to $R$.
	
	Here we should notice that if the above $\psi$ is a graded equivalence then:
	\begin{enumerate}
		\item If additionally $\phi$ is trivial then $\psi$ is a graded isomorphism.
		\item If additionally $\varphi$ is trivial then $\psi$ is an $R$-graded isometry.
		\item If additionally both $\phi$ and $\varphi$ are trivial then $\psi$ is an equivalence as Clifford algebra.
	\end{enumerate}
	
	Now, since $\psi$ is a ring homomorphism 
	\begin{equation}\label{eq:gradedequiveta1}
		\psi(u_gr)=\psi (\eta_1(g) (r)u_g)=\varphi(\eta_1 (g)(r))\lambda (\phi (g))^{-1} v_{\phi (g)}
	\end{equation}
	and on the other hand
	\begin{equation}\label{eq:gradedequiveta2}
		\psi(u_gr)=\psi (u_g)\varphi (r)=\lambda (\phi (g))^{-1} v_{\phi (g)}\varphi (r)=
		\lambda (\phi (g))^{-1}\eta_2(\phi (g))(\varphi (r))v_{\phi (g)}.
	\end{equation}
	Comparing~\eqref{eq:gradedequiveta1} and~\eqref{eq:gradedequiveta2} we get that for any $r\in R$ and any $g\in G$
	$$\varphi(\eta _1(g)(r))\lambda (\phi (g))^{-1}=
	\lambda (\phi (g))^{-1}\eta_2(\phi (g))(\varphi (r)).$$
	For simplicity, put $h=\phi (g)$ and we get
	$$\varphi(\eta _1 (\phi^{-1}(h))(r))\lambda (h)^{-1}=
	\lambda (h)^{-1}\eta_2(h)(\varphi (r))$$
	and hence by multiplying both sides of the equation on the right by $\lambda (h)$ we get that for any $h\in G$
	$$\iota _{\lambda (h)}\varphi(\eta _1 (\phi^{-1}(h)))=\eta_2(h)\varphi .$$
	Consequently,
	\begin{equation}\label{eq:actionofeta_2}
		\eta_2(h)=\iota _{\lambda (h)}\varphi(\eta _1 (\phi^{-1}(h)))\varphi ^{-1}.
	\end{equation}
	Then, by Corollary~\ref{cor:concreteactionofK} we conclude that  
	\begin{equation}\label{eq:eta_1mapstoeta_2}
		\eta _2(h)=((\lambda \varphi \phi)(\eta_1))(h).
	\end{equation}
	
	Now, for any $g_1,g_2\in G$
	\begin{equation}\label{eq:gradedequivalpha1}
		\psi (u_{g_1}u_{g_2})=\psi (\alpha_1 (g_1,g_2)u_{g_1g_2})=\varphi (\alpha_1 (g_1,g_2))\lambda(\phi (g_1g_2))^{-1}v_{\phi (g_1g_2)}
	\end{equation}
	and on the other hand
	\begin{equation}\label{eq:gradedequivalpha2}
		\begin{aligned}
			& \psi (u_{g_1}u_{g_2})=\psi (u_{g_1})\psi(u_{g_2})=
			\lambda (\phi(g_1))^{-1}v_{\phi(g_1)}\lambda (\phi(g_2))^{-1}v_{\phi(g_2)}=\\
			& =\lambda (\phi(g_1))^{-1}\eta_2 (\phi (g_1))(\lambda (\phi(g_2))^{-1})\alpha_2(\phi (g_1),\phi (g_2))v_{\phi (g_1g_2)}.
		\end{aligned}
	\end{equation}
	
	Comparing~\eqref{eq:gradedequivalpha1} and~\eqref{eq:gradedequivalpha2} we get
	$$\varphi (\alpha_1 (g_1,g_2))\lambda(\phi (g_1g_2))^{-1}=
	\lambda (\phi(g_1))^{-1}\eta_2 (\phi (g_1))(\lambda (\phi(g_2))^{-1})\alpha_2(\phi (g_1),\phi (g_2)).$$
	For simplicity put $\phi(g_1)=h_1$ and $\phi (g_2)=h_2$ then we have 
	\begin{equation}
		\alpha _2(h_1,h_2)=\eta _2(h_1)(\lambda (h_2))\lambda (h_1)\varphi(\alpha_1 (\phi^{-1}(h_1),\phi ^{-1}(h_2)))\lambda (h_1h_2)^{-1}.
	\end{equation}
	Therefore, by~\eqref{eq:actionofeta_2}, for any $h_1,h_2$
	\begin{equation}
		\begin{aligned}
			& \alpha _2(h_1,h_2)=\lambda (h_1)\varphi(\eta _1 (\phi^{-1}(h_1)))\varphi ^{-1}(\lambda (h_2))\lambda (h_1)^{-1}\lambda (h_1)\varphi(\alpha_1 (\phi^{-1}(h_1),\phi ^{-1}(h_2)))\lambda (h_1h_2)^{-1}=\\
			& =\lambda (h_1)\varphi(\eta _1 (\phi^{-1}(h_1)))\varphi ^{-1}(\lambda (h_2))\varphi(\alpha_1 (\phi^{-1}(h_1),\phi ^{-1}(h_2)))\lambda (h_1h_2)^{-1}.
		\end{aligned}
	\end{equation}
	Thus, by Corollary~\ref{cor:concreteactionofK} we conclude that 
	$$(\alpha _2,\eta _2)=(\lambda \varphi \phi)(\alpha_1,\eta_1)$$ and hence 
	$(\alpha _1,\eta _1)\in \Gamma$ and $(\alpha _2,\eta _2)\in \Gamma$ are in the same $K$ orbit.
	Moreover, by the definition of $\psi$ in~\eqref{eq:thedefinitionofPsi}
	\begin{enumerate}
		\item If $\phi$ is trivial then $\psi$ is graded isomorphism and moreover, by the above calculation 
		$(\alpha _2,\eta _2)=(\lambda \varphi )(\alpha_1,\eta_1)$, that is
		$(\alpha _1,\eta _1)\in \Gamma$ and $(\alpha _2,\eta _2)\in \Gamma$ are in the same $(R^*)^G \rtimes \text{Aut}_{\mathbb{F}}(R)$-orbit.
		\item If $\varphi$ is trivial then $\psi$ is $R$-graded isometry and moreover, by the above calculation 
		$(\alpha _2,\eta _2)=(\lambda \phi )(\alpha_1,\eta_1)$, that is
		$(\alpha _1,\eta _1)\in \Gamma$ and $(\alpha _2,\eta _2)\in \Gamma$ are in the same $(R^*)^G \rtimes \text{Aut}(G)$-orbit.
		\item If $\phi$ and $\varphi$ are both  trivial then $\psi$ is an equivalence of Clifford systems and moreover, by the above calculation 
		$(\alpha _2,\eta _2)=\lambda  (\alpha_1,\eta_1)$, that is
		$(\alpha _1,\eta _1)\in \Gamma$ and $(\alpha _2,\eta _2)\in \Gamma$ are in the same $(R^*)^G$-orbit.
	\end{enumerate}
	This concludes the ``only if" direction of Theorem A.
	
	Essentially we did also most of the work to prove the ``if" direction of Theorem A.
	
	Assume that $(\alpha _1,\eta _1)\in \Gamma$ and $(\alpha _2,\eta _2)\in \Gamma$ are in the same $K$-orbit. That is, there exist 
	$\lambda \in (R^*)^G$, $\varphi\in \text{Aut}_{\mathbb{F}}(R)$ and $\phi \in $Aut$(G)$ such that
	$$(\alpha _2,\eta _2)=(\lambda \varphi \phi)(\alpha_1,\eta_1).$$
	We define the map  $\psi :R^{\alpha_1}_{\eta_1}G\rightarrow 	R^{\alpha_2}_{\eta_2}G$ by
	\begin{equation}
		\psi \left(\sum_{g\in G}r_gu_g\right)=\sum_{g\in G}\varphi(r_g)\lambda (\phi (g))^{-1}v_{\phi (g)}.
	\end{equation}
	Then, it is clear that $\psi$ is a graded equivalence if and only if it satisfies both 
	$$\psi (u_gr)=\psi (u_g)\psi(r) 
	\text{ and } \psi (u_gu_h)=\psi (u_g)\psi (u_h)$$ for any $r\in R$, $g,h\in G$.
	By the above calculations these two conditions are satisfied.
	The ``restriction" to the finer equivalences, that is graded isomorphism, $R$-graded isometry and equivalence as Clifford systems is done similarly.

	\section{Examples}
	In this section we give several examples. In Example~\ref{example1} we are demonstrating the use of Theorem A and emphasizing the difference between the different graded relations. Then, Example~\ref{example2} and Example~\ref{example3} are demonstrating the use of Corollary B for the groups $C_4\times C_4$ and $S_3$. Lastly, in Example~\ref{example4} we show for $G\cong C_2$, $\mathbb{F}=\mathbb{Q}$ and  $L=\mathbb{Q}[\omega,\sqrt[3]{2}]$ for $\omega$ a primitive third root of unity, that the skew group rings which corresponds to the three distinct automorphisms of order $2$ in Aut$_{\mathbb{Q}}(L)$ are $\mathbb{Q}$-graded isomorphic. The notation corresponds to those in the introduction.
	\begin{example}\label{example1}
		Let $G\cong C_4\times C_4=\langle g\rangle\times \langle h\rangle$ and let $\eta$ be trivial. Consider two cohomology classes  $[\alpha],[\overline{\alpha}]\in H^2(G,\mathbb{C}^*)$ defined by
		$$[\alpha]: u_g^4=u_h^4=1,u_gu_h=iu_hu_g, \text{  and  } [\alpha] ^3=[\overline{\alpha}]: v_g^4=v_h^4=1,v_gv_h=-iv_hv_g .$$
		Then, for the twisted group algebras (here $\eta$ is trivial) $\mathbb{C}^{\alpha}G$ and $\mathbb{C}^{\bar{\alpha}}G$ it follows from Theorem A that
		\begin{enumerate}
			\item they are not $\mathbb{C}$-graded isomorphic.
			\item they are $\mathbb{C}$-graded equivalent.
			\item they are $\mathbb{R}$-graded isomorphic.
			\item they are not $\mathbb{R}$-equivalent as Clifford systems (see convention after Definition~\ref{def:Cliffordequiv}).
		\end{enumerate}
	\end{example}
	\begin{example}\label{example2}
		Let $G=C_4\times C_4$ and let $\mathbb{F}=\mathbb{R}$. For the case $L=\mathbb{R}$ we have for both $D=\mathbb{R}$ and $D=\mathbb{H}$ that the $G$-graded isomorphism classes are in one to one correspondence with $H^2(G,\mathbb{R}^*)\cong C_2\times C_2\times C_2$.
		
		For the case $L=\mathbb{C}=D$, we have two automorphisms, the trivial automorphism and the complex conjugation. Hence, the kernel $N$ of $\eta$ is either $G$ itself or a subgroup of index $2$. For $N=G$, $\eta$ is trivial and $$H^2_{\eta}(G,\mathbb{C}^*)=H^2(G,\mathbb{C}^*)\cong C_4=\langle [\alpha] \rangle$$
		where $[\alpha]$ defined in Example~\ref{example1}.
		 So, in this case (of $N=G$) there are $3$ $\mathbb{R}$-graded isomorphism classes, the group algebra $\mathbb{C}G$, the twisted group algebra which correspond to $[\alpha] ^2$, and the twisted group algebras which correspond to $[\alpha] $ and $[\bar{\alpha}]=[\alpha] ^3$ are $\mathbb{R}$-graded isomorphic. 
		
		Now, $H^2_{\eta}(G,\mathbb{C})\cong C_2\times C_2$. Therefore,  for any action $\eta$ with kernel $N$ of index $2$ there are $4$ $\mathbb{R}$-graded isomorphism classes. 
	\end{example}
	We want to give an example of non-abelian group.
	\begin{example}\label{example3}
		Let $G=S_3$ and let $\mathbb{F}=\mathbb{R}$.
		For the case $L=\mathbb{R}$ we have for both $D=\mathbb{R}$ and $D=\mathbb{H}$ that the $G$-grading classes are in one to one correspondence with $H^2(G,\mathbb{R}^*)\cong C_2$.
		
		For the case $L=\mathbb{C}=D$ we have again two automorphisms, the trivial automorphism and the complex conjugation. Hence, the kernel $N$ of $\eta$ is either $G$ itself or the unique subgroup of index $2$. For $N=G$, $\eta$ is trivial and $H^2_{\eta}(G,\mathbb{C}^*)=H^2(G,\mathbb{C}^*)$ is trivial. So, in this case there is only one $\mathbb{R}$-graded isomorphism class, the group algebra $\mathbb{C}G$.
		
		For $N$ the unique subgroup of index $2$, $H^2_{\eta}(G,\mathbb{C}^*)$ is cyclic of order $2$ where the non-trivial cohomology class is inflated from the cohomology class in Example~\ref{exintro} $(3)$. So here there are two $\mathbb{R}$-graded isomorphism classes.
	\end{example}

	\begin{example}\label{example4}
		Let $G\cong C_2=\langle \sigma \rangle$, let $\mathbb{F}=\mathbb{Q}$ and let $L=\mathbb{Q}[\omega,\sqrt[3]{2}]$ for $\omega$ a primitive third root of unity. Denote by $\eta _1,\eta _2,\eta_3$ the three distinct $\mathbb{Q}$-automorphisms of $L$ of order $2$. Then, we show that the skew group algebras 
		$L_{\eta _1}G$ and $L_{\eta _2}G$ are $\mathbb{Q}$-graded isomorphic (similarly for $L_{\eta _3}G$).
		Let $u_1,u_{\sigma}$ be an $L$-basis of $L_{\eta _1}G$ and let $v_1,v_{\sigma}$ be an $L$-basis of $L_{\eta _2}G$. Consider the map
		$$\psi : L_{\eta _1}G\rightarrow L_{\eta _2}G$$
		defined by $\psi (l_1+l_2u_{\sigma})=\eta_3 (l_1)+\eta _3(l_2)v_{\sigma}$.
		Then, noticing that $\eta _1 \eta_3=\eta_3 \eta _2$, a straightforward computation yields that $\psi$ is indeed a $\mathbb{Q}$-algebra isomorphism. 
	\end{example}

\end{document}